\RequirePackage{etex}

\documentclass[11 pt]{amsart}
\usepackage[utf8]{inputenc}

\usepackage[dvipsnames,svgnames,x11names,hyperref]{xcolor}
\definecolor{ffmblue}{HTML}{006092}

\usepackage[hyphens]{url}
\usepackage[colorlinks=true,linkcolor=Maroon,citecolor=Maroon,urlcolor=ffmblue,hypertexnames=false,linktocpage]{hyperref}
\usepackage{bookmark}
\usepackage{amsthm,thmtools,amssymb,amsmath,amscd}
\RequirePackage{etex}

\usepackage{fancyhdr}
\usepackage{esint}
\usepackage{enumerate}

\usepackage{pictexwd,dcpic}
\usepackage{graphicx}
\usepackage{caption}
\swapnumbers
\declaretheorem[name=Theorem,numberwithin=section]{thm}
\declaretheorem[name=Remark,style=remark,sibling=thm]{rem}
\declaretheorem[name=Lemma,sibling=thm]{lemma}

\declaretheorem[name=Corollary,sibling=thm]{cor}

\numberwithin{equation}{section}

\usepackage{cleveref}
\crefname{lemma}{Lemma}{Lemmata}
\crefname{prop}{Proposition}{Propositions}
\crefname{thm}{Theorem}{Theorems}
\crefname{cor}{Corollary}{Corollaries}
\crefname{defn}{Definition}{Definitions}
\crefname{example}{Example}{Examples}
\crefname{rem}{Remark}{Remarks}
\crefname{assum}{Assumption}{Assumptions}
\crefname{nota}{Notation}{Notation}
\crefname{ex}{Exercise}{Exercises}

\usepackage{autonum}

\newcommand{\ti}{\tilde}

\newcommand{\bs}{\backslash}
\newcommand{\cn}{\colon}
\newcommand{\sub}{\subset}

\newcommand{\mr}{\mathring}

\newcommand{\bbR}{\mathbb{R}}

\newcommand{\bbS}{\mathbb{S}}

\newcommand{\bbM}{\mathbb{M}}

\newcommand{\8}{\infty}

\newcommand{\al}{\alpha}
\newcommand{\be}{\beta}

\newcommand{\de}{\delta}
\newcommand{\ep}{\epsilon}
\newcommand{\ka}{\kappa}
\newcommand{\la}{\lambda}

\newcommand{\si}{\sigma}

\newcommand{\vp}{\varphi}

\newcommand{\vt}{\vartheta}
\newcommand{\Om}{\Omega}
\newcommand{\De}{\Delta}
\newcommand{\Ga}{\Gamma}
\newcommand{\Th}{\Theta}

\newcommand{\cF}{\mathcal{F}}
\newcommand{\cG}{\mathcal{G}}
\newcommand{\cH}{\mathcal{H}}
\newcommand{\cI}{\mathcal{I}}

\newcommand{\cL}{\mathcal{L}}

\newcommand{\cS}{\mathcal{S}}

\newcommand{\cU}{\mathcal{U}}

\newcommand{\del}{\partial}
\newcommand{\n}{\nabla}

\newcommand{\fa}{\forall}
\newcommand{\rt}{\sqrt}

\newcommand{\ip}[2]{\left\langle #1,#2 \right\rangle}
\newcommand{\fr}[2]{\frac{#1}{#2}}
\newcommand{\tfr}[2]{\tfrac{#1}{#2}}
\newcommand{\x}{\times}

\DeclareMathOperator{\id}{id}

\DeclareMathOperator{\osc}{osc}
\DeclareMathOperator{\const}{const}
\DeclareMathOperator{\dist}{dist}

\DeclareMathOperator{\tr}{tr}

\newcommand{\pf}[1]{\begin{proof}#1 \end{proof}}
\newcommand{\eq}[1]{\begin{equation}\begin{alignedat}{2} #1 \end{alignedat}\end{equation}}

\newcommand{\br}[1]{\left(#1\right)}
\newcommand{\abs}[1]{\lvert #1\rvert}
\newcommand{\enum}[1]{\begin{enumerate}[(i)] #1 \end{enumerate}}

\newcommand{\ra}{\rightarrow}
\newcommand{\hra}{\hookrightarrow}

\newcommand{\mt}{\mapsto}

\newcommand{\hp}{\hphantom}
\newcommand{\q}{\quad}

\begin{document}
\title{Stability from rigidity via umbilicity}
\author[J. Scheuer]{Julian Scheuer}

\address{\flushleft\parbox{\linewidth}{Goethe-Universit\"at\\ Institut f\"ur Mathematik\\ Robert-Mayer-Str.~10\\ 60325 Frankfurt\\ Germany\\ {\href{mailto:scheuer@math.uni-frankfurt.de}{scheuer@math.uni-frankfurt.de}}}}
\date{\today}
\keywords{Soap bubble theorem; Serrin's problem; Steklov problem; Stability; Alexandroff-Fenchel inequalities}
\thanks{This work was partially funded by the "Deutsche Forschungsgemeinschaft" (DFG, German research foundation); Project "Quermassintegral preserving local curvature flows"; Grant number SCHE 1879/3-1.}
\begin{abstract}	
We consider a range of geometric stability problems for hypersurfaces of spaceforms. One of the key results is an estimate relating the distance to a geodesic sphere of an embedded hypersurface with integral norms of the traceless Hessian operator of a level set function for the open set bounded by the hypersurface.
As application, we give a unified treatment of many old and new stability problems arising in geometry and analysis. Those problems ask for spherical closeness of a hypersurface, given a geometric constraint. Examples include stability in Alexandroff's soap bubble theorem in space forms, Serrin's overdetermined problem, a Steklov problem involving the bi-Laplace operator and non-convex Alexandroff-Fenchel inequalities.
\end{abstract}
\maketitle

\tableofcontents

\section{Introduction and results}

\subsection{Rigidity~vs.~Stability}

This paper is about rigidity theorems in classical hypersurface theory and their stability properties. An example is the soap bubble theorem, first proved by Alexandroff: Every constant mean curvature hypersurface embedded in Euclidean space, hyperbolic space or in a half-sphere of dimension $n+1$, $n\geq 2$, must be a round sphere, see \cite{Alexandroff:12/1962}.  Way earlier Liebmann proved the convex case \cite{Liebmann:03/1900} for hypersurfaces of the Euclidean space. Even earlier Jellet \cite{Jellet:/1853} proved the two-dimensional case for starshaped surfaces, following up the treatment of surfaces of revolution due to Delaunay \cite{Delaunay:/1841}.

A related famous result, which follows from Liebmann's theorem and the Codazzi equation, is the so-called {\it{Nabelpunktsatz}} seemingly proven by Darboux\footnote{I did not find a clean reference.}: For $n\geq 2$, every totally umbilic closed, immersed hypersurface of the $(n+1)$-dimensional Euclidean space is a round sphere.

Both of these theorems are instances of so-called {\it{rigidity results}}: A particular functional $\cF$ defined on a subclass $\cI_{N}$ of the class of closed, immersed and orientable hypersurfaces $M$ of a given Riemannian manifold $N$ dictates their global shape. Abstractly this can be described as follows. Let 
\eq{\cG = \{f\cn S\ra \bbR|~ S\sub N\}}
be the set of all functions defined on a subset of $N$. Then both of the above examples describe particular level sets of a functional
\eq{\cF\cn \cI_{N}\ra \cG.}
In Liebmann's theorem we have   
\eq{\cI_{N} = \{M\sub \bbR^{n+1}\cn M~\mbox{closed}, h>0\},\q \cF_{1}(M)=H.}
Here $h$ is the second fundamental form of $M$ and $H$ its trace with respect to the induced metric. In Alexandroff's theorem we have 
\eq{\cI_{N}=\{M\sub N\cn M~\mbox{closed, embedded}\},}
where $N$ is Euclidean, hyperbolic or spherical.
In the Nabelpunktsatz we instead use
\eq{\cF_{2}(M)=\abs{\mr{A}}^{2}=\abs{A}^{2}-\tfr{1}{n}H^{2},}
where $A$ is the Weingarten operator. Then the above theorems can be phrased as
\eq{\label{AL}\cF_{1}^{-1}(\{f=c\cn c\in \bbR\})=\{M\sub N\cn M=\mbox{sphere}\}=\cF_{2}^{-1}(\{f\equiv 0\}).}
In the course of this paper, we will encounter many more examples of such functionals.

In this abstract context, the next natural question is whether we can describe enlarged preimage sets of $\cF$. \eqref{AL} characterises the preimages of particular constant functions. We are interested in the description of sets like $\cF^{-1}(\{\|f-c\|<\ep\cn c\in \bbR\})$, for some small $\ep$ and a norm $\|\cdot\|$ to be specified. Formally, if $2^{B}$ is the power set of a set $B$, then the associated question of {\it{stability}} aims to relax the constancy in the following sense:
\begin{center}
{\it{
Is the map $\cF^{-1}\cn 2^{\cG}\ra 2^{\cI}$
continuous at all arguments, the values of which we can characterise?
}}
\end{center}
Of course, the question of continuity has to be addressed with respect to topologies to be specified within particular problems.

For example, the stability question for the Nabelpunktsatz would translate to whether a hypersurface is {\it{close}} to sphere, whenever it is {\it{almost}} umbilical. This question has received plenty of attention, initiated by early works of Reshetnyak \cite{Reshetnyak:07/1968} and Pogorelov \cite{Pogorelov:/1973} for convex surfaces with a follow up by Leichtwei\ss, see \cite{Leichtweis:08/1999}, related work by Drach \cite{Drach:/2015} and a version for hyperbolic and spherical spaces \cite{ChengZhou:07/2014}. Some optimality results are given in \cite{ChengJuarez:08/2015}. The study of almost umbilical hypersurfaces in the non-convex class seems pioneered by works of De Lellis and M\"uller \cite{De-LellisMuller:/2005, De-LellisMuller:/2006}, who give a $W^{2,2}$-closeness, resp. a $C^{0}$-closeness in terms of an $L^{2}$-pinching of $\abs{\mr{A}}^{2}$ in $\bbR^{3}$. The $W^{2,2}$-closeness result was reproved with the help of the Willmore flow in \cite{KuwertScheuer:/2021}.  Further results of almost umbilical type in terms of other $L^{p}$-norms were deduced in \cite{Perez:/2011,Roth:02/2013,Roth:/2015,RothScheuer:12/2018,Scheuer:08/2015}. 
In particular we will make use of the following estimate due to De Rosa/Gioffr\'e \cite[Thm.~1.3]{De-RosaGioffre:/2021} for a closed, connected, oriented hypersurface $M\sub\bbR^{n+1}$  and $p>n\geq 2$:
\eq{\label{RS}\dist(M,S)\leq C\|\mr{A}\|_{p},}
for some sphere $S$ with constant $C=C(n,p,\|A\|_{p},\abs{M})$. Here $\|\cdot\|_{p}$ is the $L^{p}$-norm with respect to the surface area measure.

For the soap bubble theorem most of the stability results available control the Hausdorff-distance of the hypersurface to a sphere (or an array of spheres) in terms of various pinching assumptions on the mean curvature. In the class of convex hypersurfaces such results are pioneered by Koutroufiotis \cite{Koutroufiotis:/1971} and Moore \cite{Moore:06/1973}, also see \cite{Arnold:/1993} for a rather convex geometric approach.
In the non-convex class the geometry of the stability problem is more complicated, since bubbling (closeness to a connected sum of spheres with small necks) can occur in case the curvature of the hypersurfaces is uncontrolled. In terms of a $C^{0}$-pinching condition, stability results in the Euclidean and hyperbolic space were proven by Ciraolo/Vezzoni \cite{CiraoloVezzoni:02/2018,CiraoloVezzoni:/2020}, where in both cases the authors used a quantitative analysis of Alexandroff's proof by moving planes. In \cite{CiraoloMaggi:04/2017} similar stability results were obtained, but with a proof that relies on the stability analysis of Ros' proof of the soap bubble theorem via the Heintze-Karcher inequality \cite{Ros:/1987}.

In parallel to these developments, various stability results by Magnani\-ni/Poggesi appeared relaxing the $C^{0}$-pinching condition to an $L^{1}$-pinching condition. The key is, that the soap bubble theorem can also be proved by integral methods, see \cite{Ros:/1987} and also compare \cite{Reilly:03/1982}. With the help of integral formulae, which are valid for solutions $f$ of the {\it{torsion Dirichlet problem}}
\eq{\label{Dirichlet} \bar\De f&=1\q\mbox{in}~\Om\\
					f&=0\q\mbox{on}~\del\Om,} 
Magnanini/Poggesi \cite{MagnaniniPoggesi:10/2019, MagnaniniPoggesi:01/2020,MagnaniniPoggesi:/2020} were able to prove that what they call {\it{Cauchy-Schwarz deficit}}
\eq{\abs{\mr{\bar\n}^{2}f}^{2}:=\abs{\bar\n^{2} f}^{2}-\tfr{1}{n+1}(\bar\De f)^{2}} 
can be controlled by the $L^{1}$ deficit of the mean curvature with respect to a constant. Comparison of $f$ with the squared distance function gave closeness of $\del\Om$ to a sphere after employing various estimates for harmonic functions.

The key idea of the present paper is that it is possible to control the distance to a sphere by an $L^{p}$-norm of the Cauchy-Schwarz deficit, without relying on $f$ being a solution to a PDE. This flexibility will allow us to treat a broader range of stability results, such as an extension of the results by Magnanini/Poggesi to spaceforms of constant curvature, stability results for the deficit in the Heintze-Karcher inequality in spaceforms, a double-stability result for Serrin's overdetermined problem with nonlinear source term, as well as a stability result in bi-Laplace Steklov eigenvalue problems. Everything comes with a price: The cost of this flexibility is that in particular geometric problems, some constants may be improved with methods tailor-made for the problem at hand.

Let's continue this introduction with a precise description of the results and more background information on the specific geometric problems.
 
 \subsection{Level set stability}
 A key theorem in our approach, which can be viewed as a generalisation and quantitative version of a result due to Reilly \cite[Lemma~3]{Reilly:03/1980}, that a function $f$ on a domain $\Om\sub\bbR^{n+1}$, which satisfies for some constant $L$
\eq{\bar\n^{2}f = L\bar g,\q f_{|\del\Om}=0,}
is the squared distance function to a suitable point up to normalisation constants, is as follows:

\begin{thm}\label{gen-Hessian-stable}
Let $n\geq 2$ and let $(N^{n+1},\bar g)$ be a smooth conformally flat Riemannian domain, i.e. $N\sub\bbR^{n+1}$ is open and
\eq{\bar g=e^{2\psi}\ip{\cdot}{\cdot}}
with $\psi\in C^{\8}(N)$. Let $M\sub N$ be a smooth, closed and connected hypersurface. Suppose that a one-sided open neighbourhood $\cU\Subset N$ of $M$ is foliated by level sets of a function $f\in C^{2}(\bar\cU)$, i.e.
\eq{\bar \cU = \bigcup_{0\leq t\leq \max\abs{f}}M_{t},\q M_{t} = \{\abs{f}=t\},} 
 such that $f_{|M}=0$ and $\abs{\bar\n f}_{|\bar\cU}>0$.  Let $p>n$ and suppose that for some real number $C_{0}\geq 1$ there holds
\eq{\abs{M}^{\fr{p-n}{np}}\max_{0\leq t\leq \max\abs{f}}\|A\|_{p,M_{t}}+\abs{\psi - \tfr 1n \log\abs{M}}_{0,\bar\cU} + \abs{M}^{\fr 1n}\abs{\bar\n\psi}_{0,\bar \cU}\leq C_{0}.}
Then there exists a constant $C = C(n,p,C_{0})$,
  such that
\eq{\br{\int_{\cU}\abs{\mr{\bar\n}^{2} f}^{p}}^{\fr{1}{p+1}}<\fr{1}{C}\abs{M}^{\fr{n+1-2p}{n(p+1)}}\min(\max\abs{f},\abs{M}^{\fr 1n}\min\abs{\bar\n f})^{\fr{p}{p+1}}}
implies
\eq{\dist(M,\cS)\leq \fr{C}{\abs{M}^{\fr{n-3p}{n(p+1)}}\min(\max\abs{f},\abs{M}^{\fr 1n}\min\abs{\bar\n f})^{\fr{p}{p+1}}} \br{\int_{\cU}\abs{\mr{\bar\n}^{2} f}^{p}}^{\fr{1}{p+1}}  } 
for a totally umbilical hypersurface $\cS\sub N$. 
\end{thm}

Here $\bar\n$ is the Levi-Civita connection of $\bar g$, $\abs{M}$ is the surface area of $M$, $|\cdot|_{0,\bar\cU}$ the $C^{0}$-norm in $\bar\cU$ and $\dist$ the Hausdorff distance of sets. $\cU\Subset N$ means that $\bar\cU$ is compactly contained in $N$ and $C^{k}(\bar\cU)$ is the usual space of $k$-times differentiable functions in $\cU$ where all derivatives are continuous up to the boundary.

%
%
%

\Cref{gen-Hessian-stable} will be proven with the help of an almost umbilical type theorem \cite{De-RosaGioffre:/2021} in \Cref{TS}. The key feature of \Cref{gen-Hessian-stable} is the following procedure: 

\begin{center}
{\it{Whenever a rigidity result is proven by  discarding an $L^{p}$-norm of the Cauchy-Schwarz deficit of a level set function for $\Om$, you get an associated stability result.}}
\end{center}

We will see plenty of examples throughout the following discussion and we expect further applications in the future.

\begin{rem}\label{scale-invariant}
\enum{
\item In practice, $\max\abs{f}$ and $\min\abs{\bar\n f}$ will often be under control by $C_{0}$ or other quantities deduced from the problem at hand.
\item Note that results like the one above can not hold without some dependence of the constant $C$ on the Hessian of $f$, due to the possibility of bubbling, see for example \cite[Rem.~5.2]{CiraoloVezzoni:02/2018} for further details. Hence all of the stability results in this paper will come with some sort of curvature dependence. Such dependence may in some situations be removed, when more information is available, i.e. convexity of the domain, see for example \cite{Leichtweis:08/1999} for estimates in the convex class. 
\item \label{scale-invariant-1} Note that \autoref{gen-Hessian-stable} is invariant under scaling of the metric and the function,
\eq{\bar g \mt \la \bar g,\q f\mt \mu f.}
}
\end{rem}

\subsection{The Heintze-Karcher inequality}
The following inequality, known today as the {\it{Heintze-Karcher inequality}}, was proved by Ros \cite{Ros:/1987} upon inspiration by related inequalities of Heintze/Karcher \cite{HeintzeKarcher:/1978}. It goes as follows. For a closed and mean convex hypersurface $M=\del\Om\sub \bbR^{n+1}$ we have
\eq{\int_{M}\fr{1}{H_{1}}\geq (n+1)\abs{\Om}=\int_{M}u,}
with equality precisely when $M$ is a round sphere. Here $u$ is the support function
\eq{u=\ip{x}{\nu}}
and $nH_{1} = H$.
In particular this trivially yields a new proof of Alexandroff's theorem, after invoking the formula
\eq{\abs{M}=\int_{M}uH_{1}.}
As to the soap bubble theorem, there was of course nothing new to be discovered from that method, however, Ros \cite{Ros:/1987} used this method to prove that constant $H_{k}$ hypersurfaces must be spheres. We will recall this proof in \Cref{stab} in the course of proving the stability result.
In fact, Ros proved this inequality in ambient spaces of non-negative Ricci curvature, however the equality case forces the ambient space to be Euclidean. For constant curvature spaceforms, refined versions of the Heintze-Karcher inequality have emerged during the recent years. Using flow methods, Brendle \cite{Brendle:06/2013} proved that for closed and mean-convex hypersurfaces $M=\del\Om$ of the half-sphere or the hyperbolic space there holds
\eq{\int_{M}\fr{\vt'}{H_{1}}\geq \int_{M}u,}
with equality precisely on geodesic balls. Here
\eq{u=\bar g(\vt\del_{r},\nu)}
and the metric of the ambient space is written in polar coordinates as
\eq{\bar g=dr^{2}+\vt^{2}(r)\si,}
with the round metric $\si$. This result was, among others, recovered by Qiu/Xia \cite{QiuXia:/2015} via elliptic methods, i.e. by employing a new Reilly type formula, see \eqref{HK-stable-1} for further details. Similarly as above, those inequalities can be used to reprove Alexandroff's theorem in spaceforms. But also on other occasions the Heintze-Karcher inequality has proven useful, for example for the deduction of isoperimetric type inequalities in Riemannian manifolds, see for example \cite{ScheuerXia:11/2019}.

Following the previous discussion, it becomes an interesting question whether we have an associated stability result. Namely, does almost-equality of the Heintze-Karcher inequality imply any type of closeness to the sphere? For the Euclidean case there are the recent works of Magnanini/Poggesi \cite{MagnaniniPoggesi:10/2019,MagnaniniPoggesi:01/2020}. They prove that
for a $C^{2}$-hypersurface $M=\del\Om$ there holds
\eq{\dist(M,S)\leq C\br{\int_{M}\fr{1}{H_{1}}-(n+1)\abs{\Om}}^{\fr{\tau_{n}}{2}},}
for some $\tau_{n}$ and with a constant $C$ depending on $n$, a lower inball radius bound for $\Om$ at $M$ and on the diameter of $\Om$. In \Cref{stab} we will prove the following theorem for hypersurfaces of non-Euclidean spaceforms $\bbM_{K}$ of sectional curvature $K$, where $\bbM_{1} = S^{n+1}_{+}$ is the open northern hemisphere.  The proof applies to $K=0$ too, but that result would not be new and in particular the exponents which Magnanini and Poggesi achieve in the Euclidean space are better.  In the following, if a constant depends on $\max_{M} r$, in case $K=1$ this means the constant may blow up at the equator $\{r=\pi/2\}$.

\begin{thm}\label{thm:HK}
Let $n\geq 2$, $\be >0$ and $M^{n}\sub\bbM_{K}$, $K=\pm 1$, be a closed, embedded and strictly mean convex $C^{2,\be}$-hypersurface with uniform interior sphere condition of radius $\rho$. Then there exists a constant $C$ depending on $n,\be,\abs{M}$ and on upper bounds for $\rho^{-1},\max r_{|M}$ and $\abs{A}_{\be,M}$, such that
\eq{\br{\int_{M}\fr{\vt'}{H_{1}}-\int_{M}u}<C^{-1}}
implies
\eq{\dist(M,S)\leq C\br{\int_{M}\fr{\vt'}{H_{1}}-\int_{M}u}^{\fr{1}{n+2}},}
for a suitable geodesic sphere $S$.
\end{thm}

Here $\abs{A}_{\be,M}$ is the H\"older norm of the Weingarten operator on $M$.

\subsection{Hypersurfaces of constant curvature functions}
In \cite{Alexandroff:12/1962}, Alex\-androff did not only prove the soap bubble theorem. The robustness of the moving plane method allowed to treat much more general elliptic equations, also see \cite{Voss:04/1956} and \cite{Hsiung:/1956} as well as \cite{Brendle:06/2013,DelgadinoMaggi:11/2019,KwongLeePyo:/2018,Montiel:/1999,WuXia:06/2014} fur further developments. His theorem says that every embedded hypersurface with 
\eq{\label{CFC}F(\ka_{1},\dots,\ka_{n})=\const}
is a sphere. Here $F$ is a strictly monotone function of the principal curvatures $\ka_{i}$ of the hypersurface. Particular important examples are the higher order mean curvatures
\eq{H_{k}=\fr{1}{\binom{n}{k}}\sum_{1\leq i_{1}<\dots<i_{k}\leq n}\ka_{i_{1}}\dots\ka_{i_{k}},\q H_{0}:=1,}
defined on the Garding cone
\eq{\Ga_{k}=\{\ka\in \bbR^{n}\cn 0<H_{1},\dots,0<H_{k}\},}
on which these curvature operators are elliptic.
 As the moving plane method by Alexandroff provides great flexibility in terms of the elliptic operator, it is not able to capture general immersions. The famous Wente torus \cite{Wente:/1986} is an example of a constant mean curvature closed immersed hypersurface which is not a sphere, also see \cite{HsiangTengYu:05/1983}. Consequently, it has become an interesting question, under what condition on $F$  one can obtain a soap-bubble-type theorem even for immersed hypersurfaces. It was already proved by Hopf that immersed constant mean curvature two-spheres must be embedded round spheres, but the method is very restricted to two dimensions, see \cite{Hopf:/1989}. Barbosa/do Carmo \cite{BarbosaCarmo:/1984} showed that immersed  and stable constant mean curvature hypersurfaces must be spheres and this was generalised in \cite{BarbosaCarmoEschenburg:/1988} to spaceforms and in \cite{BarbosaColares:06/1997} to $H_{k}$.
As a different approach, integral methods have proven to be equally elegant and useful. Ros \cite{Ros:/1987} reproved Alexandroff's result with $F=H_{k}$ with the help of integral methods, i.e. a Heintze-Karcher-type inequality. It was pointed out by Korevaar \cite{Korevaar:/1988}, that despite the beauty of Ros' method, his result was already covered by Alexandroff's theorem. With the help of integral methods, new rigidity results of type \eqref{CFC} for immersed hypersurfaces of spaceforms were given under various conditions, e.g. provided that two different $H_{k}$ are both constant, see \cite{Choe:05/2002} for the Euclidean case and \cite{Reilly:/1970} for the half-sphere case, or for example for hypersurfaces of constant $H_{\ell}/H_{k}$ in spaceforms \cite{FonteneleNunes:/2018,KohLee:05/2001}. 

Before moving to stability results for fully nonlinear operators, we first generalise the $L^{1}$-stability result from \cite{MagnaniniPoggesi:10/2019,MagnaniniPoggesi:01/2020} for the soap bubble theorem to spaceforms.  We prove the following:

\begin{thm}\label{thm:CMC}
Let $n\geq 2$, $\be>0$ and $M^{n}\sub \bbM_{K}$, $K=\pm1$, be a closed, embedded $C^{2,\be}$-hypersurface with uniform interior sphere condition of radius $\rho$. Then there exists a constant $C$ as in \Cref{thm:HK}, such that
\eq{\left\|\vt'\br{\fr{n\int_{M}\vt'}{(n+1)\int_{\Om}\vt'}-H_{1}}_{+}\right\|_{1,M}<C^{-1}}
implies
\eq{\dist(M,S)\leq C\left\|\vt'\br{\fr{n\int_{M}\vt'}{(n+1)\int_{\Om}\vt'}-H_{1}}_{+}\right\|_{1,M}^{\fr{1}{n+2}}}  
for a suitable geodesic sphere $S$, where $f_{+}$ denotes the positive part of a function $f$.
\end{thm}

Despite the wide variety of stability results concerning the soap bubble theorem, i.e. $F=H_{1}$, much less seems to be known when it comes to stability results for other $F$. While in the convex class there are many results available, e.g. \cite{Arnold:/1993,GroemerSchneider:/1991,Koutroufiotis:/1971,KohUm:/2001,Moore:06/1973}, 
in the non-convex class for more general $F$ there is \cite{CiraoloRoncoroniVezzoni:10/2021}, who quantified the moving plane method and obtained spherical closeness in terms of a $C^{0}$-pinching condition, generalising the corresponding results for $F=H$ from \cite{CiraoloVezzoni:02/2018,CiraoloVezzoni:/2020}. In \Cref{stab-CFC} we extend the class of stability results for non-convex hypersurfaces and prove stability for a variety of nonlinear curvature operators, employing the method of almost umbilicity. Namely we treat the case of curvature quotients. We show that the distance to a sphere can be controlled by an $L^{1}$-norm of the deviation of $F=H_{k+1}/H_{\ell}$, $0\leq \ell\leq k$, from a constant. This weakens the pinching assumption in \cite{CiraoloRoncoroniVezzoni:10/2021} from $C^{0}$ to $L^{1}$.  A particular highlight of this theorem is that in case $\ell = k$ the curvature assumptions do not imply the ellipticity of the operator $F$. This is somewhat surprising and means that there is no chance to prove this result with previous methods, such as moving planes. In the following we denote
\eq{H_{-1}=\fr{1}{H_{1}}}  and
\eq{H_{k+1,n1} = \fr{\del^{2}H_{k+1}}{\del\ka_{n}\del\ka_{1}}.}

\begin{thm}\label{thm:CFC}
Let $p>n\geq 2$, $1\leq k\leq n-1$ and $0\leq \ell\leq k$. Let $M^{n}\sub\bbM_{K}$, $K=0,1,-1$, be a closed, immersed and orientable hypersurface, such that $\ka\in \Ga_{k}$ and in case $\ell\leq k-1$ suppose $\ka\in \Ga_{k+1}$. Then there exists a constant $C$ depending on $n,p,\abs{M}$ and on upper bounds for $\max_{M}r,\abs{A}_{0,M}$ and $(\min_{M}H_{k+1,n1})^{-1}$, such that
\eq{\left\|\vt'H_{\ell-1}\br{\fr{\int_{M}\vt'H_{k}}{\int_{M}\vt'H_{\ell-1}}-\fr{H_{k+1}}{H_{\ell}}}_{+}\right\|_{1,M}<C^{-1}}
implies
\eq{\dist(M,S)\leq C\left\|\vt'H_{\ell-1}\br{\fr{\int_{M}\vt'H_{k}}{\int_{M}\vt'H_{\ell-1}}-\fr{H_{k+1}}{H_{\ell}}}_{+}\right\|_{1,M}^{\fr 1p} }
for a suitable geodesic sphere $S$.
\end{thm}

\begin{rem}
Note that for $\ell\geq 1$, this theorem quantifies the result in \cite{KohLee:05/2001}, i.e. immersions with constant quotient $H_{k+1}/H_{\ell}$ are spheres. This is because if 
\eq{\fr{H_{k+1}}{H_{\ell}}=\cH,}
then due to the Hsiung identities we have
\eq{\int_{M}\vt'H_{k}=\int_{M}uH_{k+1}=\cH\int_{M}uH_{\ell}=\cH\int_{M}\vt'H_{\ell-1}}
and hence we are comparing to the ``correct'' constant. 

Note however, that something interesting happens when $\ell=0$. In this case the correct constant must be 
\eq{H_{k+1}=\fr{\int_{M}\vt'H_{k}}{\int_{M}u}.}
Hence the theorem does provide a stability result in spaceforms for the result of Ros \cite{Ros:/1987}.
However, we do not obtain that every immersed, constant $H_{k+1}$ hypersurface, $1\leq k\leq n-2$, must be a sphere. To my knowledge this question is still open. Since this statement is wrong for $k=0$, it seems that very different methods would be required to attack this open problem. 
 \end{rem}

However we do obtain the following corollary in this direction, which says that for $k\geq 1$, a constant $H_{k+1}$ immersed closed hypersurface, on which the Heintze-Karcher inequality holds, must be a sphere.

\begin{cor}
Let $n\geq 2$ and $1\leq k\leq n-2$. Let $M^{n}\sub\bbM_{K}$, $K=0,1,-1$, be a closed, immersed and orientable hypersurface. If $H_{k+1}=\const$ and 
\eq{\int_{M}u\leq\int_{M}\fr{\vt'}{H_{1}},}
then $M$ is a round sphere.
\end{cor}

\pf{
If $H_{k+1} = \cH$, then
\eq{\cH = \fr{\int_{M}\vt'H_{k}}{\int_{M}u}\geq \fr{\int_{M}\vt'H_{k}}{\int_{M}\vt'H_{-1}}}
and hence
\eq{\br{\fr{\int_{M}\vt'H_{k}}{\int_{M}\vt'H_{-1}}-H_{k+1}}_{+} = 0.}
}

\subsection{Non-convex Alexandroff-Fenchel inequalities}

Let $\Om\Subset\bbR^{n+1}$ be a convex domain. The classical Steiner formula \cite{Steiner:/2013} says that the volume expansion of $\ep$-parallel bodies is given by
\eq{\abs{\bar\Om+\ep B}=\sum_{k=0}^{n+1}\binom{n+1}{k}W_{k}(\Om)\ep^{k}\q\fa \ep\geq 0.}
Here $B$ is the standard unit ball.
The coefficients $W_{k}(\Om)$ are called the quermassintegrals of $\bar\Om$. It can be shown easily that in case $M=\del\Om$ is of class $C^{2}$, then
\eq{W_{0}(\Om)=\abs{\Om},\q W_{k}(\Om)=\fr{1}{n+1}\int_{M}H_{k-1},\q 1\leq k\leq n+1,}
e.g. \cite{Schneider:/2014}. The classical Alexandroff-Fenchel inequalities \cite{Alexandroff:/1937b,Alexandroff:/1938} in this setting say that
\eq{\fr{W_{k+1}(\Om)}{W_{k+1}(B)}\geq \br{\fr{W_{k}(\Om)}{W_{k}(B)}}^{\fr{n-k}{n+1-k}},\q 0\leq k\leq n,}
with equality precisely if $\Om$ is a ball. For $k=0$ this is the isoperimetric inequality. With the help of inverse curvature flows, those inequalities were generalised by Guan/Li \cite{GuanLi:08/2009} to domains which are starshaped and $k$-convex, i.e. the curvatures of $M$ lie in the Garding cone $\Ga_{k}$. It is still open whether the starshapedness assumption can be dropped without replacement.

Due to the characterisation of the equality case, we may also ask for stability properties of these inequalities. In the convex case, the Hausdorff distance to a sphere in terms of the Alexandroff-Fenchel discrepancy was estimated by Schneider \cite{Schneider:01/1989} and Groemer/Schneider \cite{GroemerSchneider:/1991}, also see \cite{Ivaki:07/2015,Ivaki:05/2016} for similar stability results in the class of convex bodies.  However, except for the case of the isoperimetric inequality, for which the stability question has been studied extensively, e.g. \cite{FuscoMaggiPratelli:/2008,FigalliMaggiPratelli:10/2010,Ivaki:10/2014}, I am not aware of a stability result for the non-convex Alexandroff-Fenchel inequalities of Guan/Li. In \Cref{stab-AF} we will combine a curvature flow with the almost umbilicity theorem from \cite{De-RosaGioffre:/2021} to prove a stability result for the non-convex Alexandroff-Fenchel inequalities. 

\begin{thm}\label{thm:AF}
Let $n\geq 2$ and $1\leq k\leq n$. Let $M^{n}\sub\bbM_{0}$ be a closed and starshaped $C^{2}$-hypersurface, such that $\ka\in \Ga_{k}$. Then there exists a constant $C$, depending on $n$, $\abs{M}$ and on upper bounds for $(\min_{M}r)^{-1},\abs{r}_{2,M}$ and $(\min_{M}\dist(\ka,\del\Ga_{k}))^{-1}$, such that
\eq{\fr{W_{k+1}(\Om)}{W_{k+1}(B)}-\br{\fr{W_{k}(\Om)}{W_{k}(B)}}^{\fr{n-k}{n-k+1}}<C^{-1}}
 implies
\eq{\dist(M,S)\leq C\br{\fr{W_{k+1}(\Om)}{W_{k+1}(B)}-\br{\fr{W_{k}(\Om)}{W_{k}(B)}}^{\fr{n-k}{n-k+1}}}^{\fr{1}{2(n+1)}} }
for a sphere $S$.
\end{thm}

\subsection{Serrin's overdetermined problem}
Very related to the soap bubble theorem and the Heintze-Karcher inequality are overdetermined elliptic problems. From the unique solvability of the Dirichlet problem
\eq{\label{Diri}\bar\De f&=1\q\mbox{in}~\Om\\
		f&=0\q\mbox{on}~\del\Om}
it is a natural question to ask how much information a further condition, such as a Neumann condition, gives. Interestingly, Serrin \cite{Serrin:01/1971} proved that the additional condition
\eq{\del_{\nu}f=\const}
implies that $f$ is, up to some normalisation, the squared radial distance and that $\Om$ must be a ball. Immediately after this result appeared, Weinberger \cite{Weinberger:01/1971} gave a simplified proof by showing that a good choice of test function yield the Cauchy-Schwarz deficit of $f$ to vanish, from which the result follows. The stability results of the present paper yields into a similar direction, as we use a quantitative version of this argument employing \Cref{gen-Hessian-stable}. Serrin's result was generalised to the hyperbolic space and the hemisphere by Kumaresan/Prajapat \cite{KumaresanPrajapat:/1998} via moving planes. Using arguments like those of Weinberger, similar results in a class of warped products and in spaceforms were obtained in \cite{QiuXia:/2017}, \cite{Roncoroni:06/2018} and in \cite{CiraoloVezzoni:07/2019}.

The stability question associated to Serrin's problem asks, whether 
\eq{\cF(M)=\abs{\del_{\nu}f-R}}
can serve as a pinching quantity for spherical closeness. In the Euclidean space, under a $C^{0}$-condition on $\abs{\bar\n f}$ this was answered affirmatively in \cite{AftalionBuscaReichel:/1999,BrandoliniNitschSalaniTrombetti:09/2008,CiraoloMagnaniniVespri:08/2016}. When it comes to pinching in terms of an integral quantity Magnanini/Poggesi \cite{MagnaniniPoggesi:/2020} have shown the following identity for 
solutions of \eqref{Diri}:
\eq{\label{MaPo}\int_{\Om}(-f)\abs{\mr{\bar\n}^{2}f}^{2}=\fr 12\int_{M}((\del_{\nu}f)^{2}-R^{2})(\del_{\nu}f-\del_{\nu}q),}
where $q$ is any squared distance function and $R$ is a suitable constant. A similar approach was taken by Feldman \cite{Feldman:/2018} and an interesting alternative approach via quermassintegrals is in \cite{BrandoliniNitschSalaniTrombetti:08/2008}. This beautiful formula shows that $\cF$ pinches the Cauchy-Schwarz deficit and hence we could directly apply \Cref{gen-Hessian-stable} to obtain spherical closeness, although the latter is not what is done in \cite{MagnaniniPoggesi:/2020}, but the authors proceed via harmonic functions; note that $h=f-q$ is harmonic. Although their approach is well suitable for Serrin's problem, the approach via \Cref{gen-Hessian-stable} seems to be more flexible when it comes to other overdetermined problems. In particular, the flexibility of \Cref{gen-Hessian-stable} allows us to prove the following {\it{double stability result}} for Serrin's problem, which allows for nonlinear right hand sides.

\begin{thm}\label{thm:Serrin}
Let $n\geq 2$, $\Om\Subset\bbM_{0}$ be an open set with connected $C^{2}$-boundary, uniform interior sphere condition with radius $\rho$ and $\abs{\del\Om}=1$. Let $\phi\in C^{2}((-\8,0])$ be positive and suppose that $f\in C^{2}(\bar\Om)$ satisfies
\eq{\bar\De f &= \phi(f)\q&&\mbox{in}~\Om\\
		f&=0\q&&\mbox{on}~\del\Om.}
There exists a constant $C$ depending on $n$ and upper bounds for the quantities $\abs{\bar\n^{2}f}_{0,\Om},\abs{\phi}_{0,\Om},\abs{1/\phi}_{0,\Om}$ and $\rho^{-1}$, and a round sphere $S$, such that
\eq{\|\abs{\bar\n f} - R\|_{1,\del\Om}+\|\phi-\phi(0)\|_{1,\Om}+\|\Phi-f\phi\|_{1,\Om}<C^{-1}}
implies
\eq{\label{thm:Serrin-1}\dist(\del\Om,S)\leq C\br{\|\abs{\bar\n f} - R\|_{1,\del\Om}+\|\phi-\phi(0)\|_{1,\Om}+\|\Phi-f\phi\|_{1,\Om}}^{\fr{1}{n+2}},}
where
\eq{R = \fr{1}{\abs{M}}\int_{\Om}\phi(f), \q \Phi(f) = -\int_{f}^{0}\phi(s)~ds.}
\end{thm}

\begin{rem}
\enum{
\item Of course, the second derivative of $f$ can be estimated in terms of $\phi$ and $\Om$ using elliptic estimates. But as we do not aim to keep track on how $\Om$ enters precisely in this estimate, we make the statement as general as possible.
\item One and a half years after the first version of this paper appeared, Chao Xia and JS proved a version of this stability result in spaceforms with less explicit exponents, \cite{ScheuerXia:10/2022}.
}
\end{rem}

\subsection{Estimates for Steklov eigenvalues}
The Steklov eigenvalue problem, introduced by Vladimir Steklov \cite{Steklov:/1902} models the steady state temperature in a domain $\Om\sub\bbR^{n+1}$, such that the flux at the boundary is proportional to the temperature, i.e.
\eq{\bar\De f &= 0\q&&\mbox{in}~\Om\\
		\del_{\nu}f &= \mu f\q&&\mbox{on}~M=\del\Om.}
The connection of the Steklov problem to free boundary minimal surfaces was treated by Fraser/Schoen \cite{FraserSchoen:03/2011,FraserSchoen:03/2016}. Various upper bounds for the first Steklov eigenvalue $\mu_{1}$ have been obtained in the literature, for example \cite{IliasMakhoul:10/2011,Roth:/2020,WangXia:10/2009}.  For example we have the result of Wang/Xia \cite{WangXia:10/2009}, which says that a uniformly convex domain with a lower bound $c$ on the principal curvatures of $M$ yields the estimate
\eq{\label{Steklov}\mu_{1}\leq \fr{\rt{\la_{1}(M)}}{nc}\br{\rt{\la_{1}(M)}+\rt{\la_{1}(M)-nc^{2}}},}
with equality precisely if $M$ is a round sphere. Here $\la_{1}$ is the first nonzero eigenvalue of the Laplacian on $M$. We refer to the survey \cite{GirouardPolterovich:/2017} for a much broader description of this vastly researched topic. Although it is possible to obtain a stability result for \eqref{Steklov} with the help of an almost-CMC result, here we will focus on another problem, that can be solved by direct application of \Cref{gen-Hessian-stable}. 
Namely we prove a stability result for a Steklov problem involving the bi-Laplace:
\eq{\label{BiLaplace}\bar\De^{2}w & = 0\q\mbox{in}~\Om\\
		w=\bar\De w - \mu_{1}\del_{\nu}w&=0\q\mbox{on}~\del\Om.}
This problem arises from applications in elasticity.
The following theorem is a stability version of the rigidity result \cite[Thm.~1.2]{WangXia:10/2009}, that whenever $H_{1}\geq c$ and $\mu_{1}=(n+1)c$, then $\del\Om$ is a round sphere. Note that the following theorem even improves this rigidity result, in the sense that we do not need a priori to assume the lower bound on $H_{1}$.

\begin{thm}\label{thm:Steklov}
Let $n\geq 2$ and  $\Om\Subset \bbM_{0}$ be an open set with smooth and connected boundary and $\abs{\del\Om}=1$. Let $\mu_{1}$ be the first nonzero eigenvalue of the problem \eqref{BiLaplace} and suppose $w$ is an eigenfunction normalised to
\eq{\min_{\del\Om}\abs{\bar\n w}=1.}
Then there exists a constant $C$ depending on $n$ and on an upper bound for $\abs{w}_{2,0,\Om}$, such that 
\eq{\|(\mu_{1}-(n+1)H_{1})_{+}\|_{1,\del\Om}<C^{-1}}
implies
\eq{\dist(\del\Om,S)\leq C\|(\mu_{1}-(n+1)H_{1})_{+}\|_{1,\del\Om}^{\fr{1}{n+2}} }
for some round sphere $S$.

\end{thm}

\section{Definitions and Notation}\label{sec:notation}
Let us collect some notation we use throughout the paper. 
Let $N\sub \bbR^{n+1}$ be open and $\bar g$ be a Riemannian metric on $N$. Then we call the pair $(N,\bar g)$ a Riemannian domain. Let $\bar\n$ denote the Levi-Civita connection of $\bar g$.

For a hypersurface $x\cn M\hra N$, induced geometric quantities such as the induced metric and its connection are distinguished from their ambient counterparts by missing an overbar, e.g. $g=\bar g^{*}$, $\n$ and $\De$.
We write $\abs{M}$ for the volume of the Riemannian manifold $(M,g)$.
The second fundamental form $h$ of $M$ in $N$ is defined by the Gaussian formula
\eq{\bar\n _{X}Y = \n_{X}Y - h(X,Y)\nu,}
where $\nu$ is a given local smooth normal, which will be chosen as ``outward pointing'', whenever such a choice is possible. For the associated Weingarten operator we write $A$, i.e. there holds
\eq{h(X,Y) = g(A(X),Y) = g(X,A(Y)).}
This operator is $g$-selfadjoint and its eigenvalues are called the {\it{principal curvatures}} usually written as 
\eq{\ka = (\ka_{1},\dots,\ka_{n}).}
A key quantity in this paper is the {\it{traceless Weingarten operator}}
\eq{\mr{A} = A - \fr{1}{n}H \id,}
where $H = \tr A = g^{ij}h_{ij}$ is the mean curvature, indices denote components with respect to a local frame
$(e_{i})_{1\leq i\leq n}$,
\eq{g_{ij} = g(e_{i},e_{j}),\q h_{ij} = h(e_{i},e_{j}),}
and where $(g^{ij})$ is the ``inverse'' of $g$, i.e.
\eq{\de^{i}_{j} = g^{ik}g_{kj}}
with the Kronecker-delta $\de^{i}_{j}$.
$M$ is called {\it{totally umbilic}}, if $\mr{A} = 0$. It is well known that in spaceforms all closed, totally umbilic hypersurfaces are geodesic spheres. Let us also note that latin indices denote components for hypersurfaces $M\sub N$, and we use greek indices to distinguish indices of ambient quantities defined on $N$, e.g
\eq{\bar{g}_{\al\be} = \bar g(e_{\al},e_{\be}).}

For a set $U\sub N$ and $T$ a tensor field defined on $U$  we use the standard essential sup-norm
\eq{\|T\|_{\8,U} = \sup_{U}~\abs{T},}
where norms of tensors are formed with respect to the metric at hand.  
For a tensor $T$ on a submanifold $M\sub N$ and $1\leq p<\8$ we define the $L^{p}$-norm
\eq{\|T\|_{p,M} = \br{\int_{M}\abs{T}^{p}}^{\fr 1p}. }
We avoid to write the volume element in such integrals, as the domain is always understood to be equipped with the standard volume element coming from the induced metric in question. For $\be>0$, we denote by $\abs{T}_{\be,M}$ the H\"older norm with exponent $\be$, which can simply be defined using coordinate charts and the standard definition of H\"older norms in Euclidean space. for $\be=0$ we define $\abs{T}_{0,M}$ to be the standard sup-norm. Similarly, $\abs{T}_{k,\be,M}$ is the H\"older norm of $C^{k,\be}$. 

Often the Riemannian domains we are working with are models for the $(n+1)$-dimensional Euclidean space, the hyperbolic space and the half-sphere, i.e. the upper hemisphere of the standard round unit sphere of the Euclidean space. We denote all of them by $\bbM_{K}$, where $K$ is the constant sectional curvature of the respective space. Introducing polar coordinates around a respective origin in these spaces, the ambient metrics $\bar g$ can all be written as
\eq{\bar g=dr^{2}+\vt^{2}(r)\si,
} 
where $\si$ is the round metric on the sphere $\bbS^{n}$ and 
\eq{\vt(r)=\begin{cases} r,& K=0\\
					\sin r,& K=1\\
					\sinh r, &K=-1.
\end{cases}
	}
	
For an immersed, orientable hypersurface $M\sub \bbM_{K}$ the higher order mean curvatures
\eq{H_{k}=\fr{1}{\binom{n}{k}}\sum_{1\leq i_{1}<\dots<i_{k}\leq n}\ka_{i_{1}}\dots\ka_{i_{k}}}
are elliptic operators on the Garding cones
\eq{\Ga_{k}=\{\ka\in \bbR^{n}\cn 0<H_{1},\dots,0<H_{k}\}.}
Instead of viewing those operators as acting on the principal curvatures, they can be viewed as functions $F$ of the Weingarten operator $(h^{i}_{j})$ and for such we denote
\eq{F^{i}_{j}=\fr{\del F}{\del h^{j}_{i}}.}
If $F=H_{k}$, $1\leq k\leq n$, then the tensor $(F^{i}_{j})$ is divergence free, see \cite[Lemma~5.8]{Gerhardt:/2007}. 
Denoting by $\Th$ a primitive of $\vt$, $\Th'=\vt$, the following Hsiung type formula, compare \cite{Hsiung:/1954}, follow from integration of $F^{ij}\n^{2}_{ij}\Theta$ along a hypersurface $M\sub\bbM_{K}$:
\eq{\label{Hsiung}\int_{M}\vt'H_{k}=\int_{M}uH_{k+1},\q k=0,\dots,n-1.}
Here $u$ is the support function
\eq{u=\bar{g}(\vt\del_{r},\nu)}
with a smooth unit normal $\nu$ and we use the convention $H_{0}=1$.

\section{Proof of almost umbilicity estimates for level sets}\label{TS}

In this section we prove \cref{gen-Hessian-stable}.

\pf{
%

Let $\cU$ be given as stated. Due to \autoref{scale-invariant} (\ref{scale-invariant-1}) without loss of generality we may assume
 \eq{\abs{M} = 1.}
Without loss of generality we may assume $f_{|\cU}>0$, for otherwise consider $-f$ instead.
Since all $M_{t} = \{f=t\}$ are $C^{2}$-level sets of $f$, we can relate the Hessian of $f$ to the second fundamental form of $M_{t}$. A similar calculation appeared in \cite{ColdingMinicozzi:03/2014}: Let $(x_{i})_{1\leq i\leq n}$ be local coordinate vectors for $M_{t}$. The unit normal $\nu$, which is pointing out of $\cU$, is
\eq{\nu=-\fr{\bar\n f}{\abs{\bar\n f}}.}
 Differentiate $f$ along these coordinates and obtain
 \eq{0=\bar\n f(x_{i}),\q 0=\bar\n^{2}f(x_{i},x_{j})+\bar\n f(x_{ij})=\bar\n^{2}f(x_{i},x_{j})-\bar\n f(\nu)h_{ij},}
and hence
\eq{\label{gen-Hessian-stable-6}\bar\n^{2}f(x_{i},x_{j})=-\abs{\bar\n f}h_{ij}.}
We obtain
\eq{-\abs{\bar\n f}\mr{h}_{ij}&=-\abs{\bar\n f}\br{h_{ij}-\tfr{1}{n}H g_{ij}}\\
						&=\bar\n^{2}f(x_{i},x_{j})-\tfr{1}{n}g^{kl}\bar\n^{2}f(x_{k},x_{l})g_{ij}\\
						&=\mr{\bar\n}^{2}f(x_{i},x_{j})+\tfr{1}{n+1}\bar\De f \bar g(x_{i},x_{j})\\
						&\hp{=}-\tfr{1}{n}g^{kl}\br{\mr{\bar\n}^{2}f(x_{k},x_{l})+\tfr{1}{n+1}\bar\De f\bar g(x_{k},x_{l})}g_{ij}\\
						&=\mr{\bar\n}^{2}f(x_{i},x_{j})+\tfr{1}{n+1}\bar\De f \bar g(x_{i},x_{j})-\tfr{1}{n}g^{kl}\mr{\bar\n}^{2}f(x_{k},x_{l})g_{ij}-\tfr{1}{n+1}\bar\De f g_{ij}\\
						&=\mr{\bar\n}^{2}f(x_{i},x_{j})+\tfr{1}{n}\mr{\bar\n}^{2}f(\nu,\nu)g_{ij},}
where we used
\eq{g^{kl}x^{\al}_{k}x^{\be}_{l} = \bar g^{\al\be} - \nu^{\al}\nu^{\be}.}
In turn we get the pointwise relation
\eq{\label{gen-Hessian-stable-2}\abs{\bar\n f}^{2}\abs{\mr{A}}^{2}=\mr{\bar\n}^{2}f(x_{i},x_{j})\mr{\bar\n}^{2}f(x_{k},x_{l})g^{ik}g^{jl}-\tfr{1}{n}(\mr{\bar\n}^{2}f(\nu,\nu))^{2}.  }
The level sets can be constructed as the flow hypersurfaces of the flow
\eq{\del_{t}\Phi(t,\xi)&=\fr{\bar\n f(\Phi(t,\xi))}{\abs{\bar\n f(\Phi(t,\xi))}^{2}}\\
 				\Phi(0,\xi)&=\xi,}
which makes $\cU$ diffeomorphic to $(0,\max_{\bar\cU}f)\x M$.
We have to make sure that the surface areas of the level sets $M_{t}$ are comparable to $\abs{M}=1$. To accomplish this, note that the evolution of the area element under a hypersurface variation with speed $\del_{t}\Phi$ is given by
\eq{\label{pf:gen-Hessian-stable-9}\del_{t}\rt{\det g}=H\bar g(\del_{t}\Phi,\nu)\rt{\det g}=-\fr{H}{\abs{\bar\n f(\Phi(t,\cdot))}}\rt{\det g}.}
Hence we have, with $q=p/(p-1)$, 
\eq{\del_{t}\abs{M_{t}} \leq \fr{\rt{n}}{\min_{\bar\cU}\abs{\bar\n f}}\int_{M_{t}}\abs{A}&\leq \fr{\rt{n}}{\min_{\bar\cU}\abs{\bar\n f}}\br{\int_{M_{t}}\abs{A}^{p}}^{\fr 1p}\abs{M_{t}}^{\fr 1q} \\
				&\leq \fr{\rt{n}C_{0}}{\min_{\bar\cU}\abs{\bar\n f}}\max(1,\abs{M_{t}}) }
and a similar estimate from below. We obtain
\eq{1-\fr{\rt{n}C_{0}}{\min_{\bar\cU}\abs{\bar\n f}}t\leq\abs{M_{t}}\leq e^{\fr{\rt{n}C_{0}}{\min_{\bar\cU}\abs{\bar\n f}}t}.}
Define
\eq{T_{1} = \min\br{\tfr{\min_{\bar\cU}\abs{\bar\n f}}{2\rt{n}C_{0}},\max_{\bar\cU}f}\geq  \fr{1}{2\rt{n}C_{0}}\min\br{\min_{\bar\cU}\abs{\bar\n f},\max_{\bar\cU}f},}
then 
\eq{\tfr 12\leq \abs{M_{t}}\leq 2\q\fa 0<t<T_{1}.}

For all $t_{0}\in (0,T_{1})$, let us integrate \eqref{gen-Hessian-stable-2} and use the co-area formula,
\eq{\int_{0}^{t_{0}}\int_{M_{s}} \abs{\bar\n f}^{p-1}\abs{\mr{A}}^{p}~ ds&\leq \int_{0}^{t_{0}}\int_{M_{s}}\fr{\abs{\mr{\bar\n}^{2}f}^{p}}{\abs{\bar\n f}} ds.}
Hence
\eq{\cL^{1}\br{\left\{s\in (0,t_{0})\cn\int_{M_{s}}\abs{\mr A}^{p}>\fr{2}{t_{0}\min_{\bar\cU}\abs{\bar\n f}^{p-1}}\int_{\cU}\abs{\mr{\bar\n}^{2}f}^{p}\right\}}\leq \fr{t_{0}}{2},}
where $\cL^{1}$ is the one-dimensional Lebesgue measure.
Define
\eq{t_{0} := 2\min(\min_{\bar\cU}\abs{\bar\n f},\max_{\bar\cU} f)^{\fr{1}{p+1}}\|\mr{\bar\n}^{2}f\|_{p,\cU}^{\fr{p}{p+1}},}
while we have to assume that this is less than $T_{1}$ and hence we demand, as a first condition,
\eq{\|\mr{\bar\n}^{2}f\|_{p,\cU}^{\fr{p}{p+1}}\leq \fr{1}{4\rt{n}C_{0}}\min\br{\min_{\bar\cU}\abs{\bar\n f},\max_{\bar\cU}f}^{\fr{p}{p+1}}.}
Thus there exists $s\in (0,t_{0})$, such that $M_{s}$ has the property
\eq{\label{pf:stability-1}\|\mr{A}\|_{p,M_{s}}	\leq \fr{1}{\min\br{\min_{\bar\cU}\abs{\bar\n f},\max_{\bar\cU}f}^{\fr{p}{p+1}}}\|\mr{\bar\n}^{2}f\|_{p,\cU}^{\fr{p}{p+1}}.}

We want to apply \cite[Thm.~1.3]{De-RosaGioffre:/2021}, which applies to hypersurfaces of the Euclidean space. Therefore, for a better distinction, let us denote by $\hat M_{s}$ the hypersurface $M_{s}$ equipped with the metric induced from the Euclidean metric, i.e.
\eq{\hat g = e^{-2\psi}g, }
and also the other geometric quantities on $\hat M_{s}$ will be wearing a hat.  For the Weingarten operator we have
\eq{e^{\psi}A = \hat A + e^{\psi}d\psi(\nu)\id,}
see \cite[Prop.~1.1.11]{Gerhardt:/2006} 
and hence we obtain
\eq{\abs{\mr{\hat A}}^{p}\rt{\det \hat g} = e^{(p-n)\psi }\abs{\mr{A}}^{p}\rt{\det g},}
as well as
\eq{\|\hat A\|_{p,\hat M_{s}}=\br{\int_{ M_{s}}e^{(p-n)\psi}\abs{A - d\psi(\nu)\id}^{p}}^{\fr 1p}\leq C,}
where $C$ depends on $n$, $C_{0}$ and $\abs{\psi}_{1,\bar\cU}$.
In addition we have
\eq{\|\mr{\hat A}\|_{p,\hat M_{s}}\leq  \fr{C}{\min\br{\min_{\bar\cU}\abs{\bar\n f},\max_{\bar\cU}f}^{\fr{p}{p+1}}}\|\mr{\bar\n}^{2}f\|_{p,\cU}^{\fr{p}{p+1}}.}
We have to scale to unit sphere area, i.e. we define
\eq{\ti M_{s} = \la\hat M_{s},}
where, with $C = C(n,\abs{\psi}_{0,\bar \cU})$, 
\eq{2^{-\fr 1n}C^{-1}\leq\fr{C^{-1}}{\abs{M_{s}}^{ 1/n}}\leq\la = \br{\fr{\abs{\bbS^{n}}}{\abs{\hat M_{s}}}}^{\fr 1n} \leq \fr{C}{\abs{M_{s}}^{1/n}}\leq 2^{\fr 1n}C.}
Then,
\eq{\|\mr{\ti A}\|_{p,\ti M_{s}}=\la^{\fr{n-p}{p}}\|\mr{\hat A}\|_{p,\hat M_{s}}. } According to \cite[Thm.~1.3]{De-RosaGioffre:/2021}, there exist numbers $\de_{1}$ and $C_{1}$, depending only on $n$, $p$ and $C_{0}$, such that whenever
\eq{\|\mr{\ti A}\|_{p,\ti M_{s}}\leq \de_{1},}
then $\ti M_{s}$ is a graph over the unit sphere centred at (without loss of generality) the origin and the graph function $\ti\vp$ (i.e. the distance to the origin) satisfies
\eq{\|\log\ti\vp\|_{2,p,\bbS^{n}}\leq C_{1}\|\mr{\ti A}\|_{p,\ti M_{s}}.}
The threshold $\de_{1}$ can be achieved by simply using the smallness assumption on $\|\mr{\bar\n}^{2}f\|_{p,\cU}$.
Scaling back to $\hat M_{s}\sub\cU$, we obtain with $\hat \vp = \la^{-1}\ti\vp$
\eq{\label{pf:stability-2}\|\log\la + \log\hat \vp\|_{2,p,\bbS^{n}}\leq C\la^{\fr{n-p}{p}}\|\mr{ \hat A}\|_{p,\hat  M_{s}}.}
In particular we obtain for the Hausdorff distance,
\eq{\dist(\hat M_{s},S_{\la^{-1}})\leq C\la^{\fr{n-p}{p}}\|\mr{\hat A}\|_{p,\hat M_{s}}\leq \fr{C}{\min\br{\min_{\bar\cU}\abs{\bar\n f},\max_{\bar\cU}f}^{\fr{p}{p+1}}}\|\mr{\bar\n}^{2}f\|_{p,\cU}^{\fr{p}{p+1}} }
and we obtain the same estimate for the Hausdorff distance of $M_{s}$ to the conformal image $\cS$ of $S_{\la^{-1}}$, due to the equivalence of the metrics on $\bar\cU$.
It remains to transfer this estimate to $M$. Recall that $s\leq t_{0}$
 and that every point $\xi\in M$ is connected to a unique point in $M_{s}$ by the curve $\Phi(\cdot,\xi)$,
where
\eq{\xi=\Phi(0,\xi)\in M,\q \Phi(s,\xi)\in M_{s}.}
Along this curve, we estimate $\Phi$ by
\eq{\label{pf:gen-Hessian-stable-8}d_{\bar g}(\Phi(s,\xi),\Phi(0,\xi))\leq \fr{Ct_{0}}{\min(\min_{\bar\cU}\abs{\bar\n f},\max_{\bar\cU}f)}}
and hence
\eq{\dist(M,\cS)\leq \fr{Ct_{0}}{\min(\min_{\bar\cU}\abs{\bar\n f},\max_{\bar\cU}f)},}
which finishes the proof.
}

\begin{rem}
With more regularity, it would be possible to obtain more information on the zero-set $M$. For example, under a Hölder continuity assumption on the gradient of $f$, one can deduce that $M$ is starshaped. The available regularity depends on the concrete geometric problem. In \cite{ScheuerXia:10/2022}, the method of this proof was used in particular geometric settings.
\end{rem}

\section{Proofs of the geometric stability results}\label{stab}

\subsection{Stability in the Heintze-Karcher inequality}

The proof of \Cref{thm:HK} is built upon solving

\eq{\label{HK-Diri}\bar\De f+(n+1)Kf&=1\q\mbox{in}~\Om\\
					f&=0\q\mbox{on}~\del\Om}
and applying \Cref{gen-Hessian-stable}. Therefore we need a lower bound on $\abs{\bar\n f}$ along $\del\Om$. For this purpose we prove the following quantified Hopf boundary lemma in some greater generality for further applications later on.

\begin{lemma}\label{quant-Hopf}
Let $\Om\Subset \bbM_{K}$, $K=0,1,-1$, be an open set satisfying a uniform interior sphere condition with radius $\rho$. Let $f$ satisfy 
\eq{\bar\De f+af&\geq b\q\mbox{in}~\Om\\
					f&\leq 0\q\mbox{on}~\del\Om,}
where $b>0$ and $a\in C^{0}(\bar\Om)$ satisfies
\eq{a\leq (n+1)\max(0,K).} 
Then $f<0$ in $\Om$ and for all $x_{0}\in \del\Om$ with $f(x_{0})=0$ we have
\eq{\abs{\bar\n f(x_{0})}\geq \ep_{0}>0, }
where $\ep_{0}$ depends on lower bounds of $n^{-1},b,\rho$, and $(\max_{\bar\Om}\abs{a})^{-1}$. 
\end{lemma}

\pf{
The proof follows well known comparison strategies, e.g. as in \cite[Lemma~3.4]{GilbargTrudinger:/2001} in conjunction with the choice of a good test function.

(i) First we prove that $f<0$ in $\Om$. For $K\leq 0$ this follows from the strong maximum principle applied to $f$. If $K= 1$, we have
\eq{\bar\De\vt' = -(n+1)\vt'.}
By assumption $\bar\Om$ lies in the northern open hemisphere and hence $\vt'>0$. 

The function $z=f/\vt'$ solves
\eq{\bar\De z &= \fr{\bar\De f}{\vt'}-\fr{f}{\vt'^{2}}\bar\De \vt'-\fr{2}{\vt'}\ip{\bar\n z}{\bar\n \vt'}\\
			&\geq\fr{b-af}{\vt'}+(n+1)\fr{f}{\vt'}-\fr{2}{\vt'}\ip{\bar\n z}{\bar\n \vt'}\\
			&= \fr{b}{\vt'}-\fr{2}{\vt'}\ip{\bar\n z}{\bar\n\vt'} + ((n+1)-a)z.}
Since $z\leq 0$ on $M$, $z$ and hence $f$ are both negative in $\Om$ due to the strong maximum principle.

(ii) Let $B=B_{\rho}(y_{0})\sub \Om$, $0<\rho<\pi/4$, be an interior ball touching at $x_{0}\in\del\Om$. By an ambient isometry we may assume that our polar coordinates are centred at $y_{0}$ and in these coordinates the metric is given by
\eq{\bar g=dr^{2}+\vt^{2}(r)\si.}
 Let $\Th$ be a primitive of $\vt$ vanishing on $\del B$, i.e.
 \eq{\Th(r) = \int_{0}^{r}\vt(s)~ds - \int_{0}^{\rho}\vt(s)~ds .}
  $\Th$ satisfies
 \eq{\bar\De\Th=(n+1)\vt'.}
For $\de>0$ the function 
 \eq{w=f-\de\Th}
 satisfies
 \eq{\bar\De w=\bar\De f - \de(n+1)\vt'&\geq b-af - \de(n+1)\vt'\\
 					&=b-aw-\de a\Th-\de(n+1)\vt'\\
					&\geq \tfr b2 - aw, }
if $\de=\de(n^{-1},b,(\max_{\bar\Om}\abs{a})^{-1})$ is chosen suitably small. As $w=f\leq 0$ on $\del B$, we may apply part (i) of this proof to $w$ and obtain $w\leq 0$ in $B$. Hence
\eq{f_{|B}\leq \de\Th_{|B}.}
Now define
\eq{v(x)=e^{-C \Th(r)}-1,}
where $r=r(x)$.
Then in the region $\rho/2<r<\rho$ we can estimate
\eq{\bar\De v&=-C e^{-C \Th}\bar\De\Th+C^{2}e^{-C \Th}\abs{\bar\n\Th}^{2}\\
		&=e^{-C\Th}(-C (n+1)\vt'+C^{2}\abs{\bar\n\Th}^{2})\\
		&\geq e^{-C\Th}(-C (n+1)\vt'+C^{2}\vt^{2}(\tfr{\rho}{2}))\\
		&\geq 0,}
provided $C=C(n,\rho^{-1})$ is chosen large enough. Since
\eq{f_{|\del B_{\rho/2}}\leq -\de\int_{\rho/2}^{\rho}\vt(s)~ds,}
there exists a constant $\ep = \ep(\de,\rho,b,(\max \abs{a})^{-1})$, such that
\eq{f_{|\del B_{\rho/2}}+\ep v_{|\del B_{\rho/2}}\leq 0\q\mbox{and}\q \ep a v_{|B}\geq -\tfr b2.}
Hence we have 
\eq{\bar\De(f+\ep v)+ a(f+\ep v)\geq b + \ep\bar\De v + \ep av \geq \tfr b2.}
Step (i) applied to the domain $\bar B_{\rho}\bs B_{\rho/2}$ gives
\eq{f+\ep v\leq 0\q\mbox{in}~\bar B_{\rho}\bs B_{\rho/2}. }
Taking the outward pointing derivative at $x_{0}$ we get
\eq{0\leq \del_{\nu}f+\ep\del_{\nu}v=\del_{\nu}f-\ep C e^{-C\Th(\rho)}\vt(\rho)}
and
\eq{\del_{\nu}f\geq \ep C\vt(\rho).}
}

\begin{proof}[Proof of \Cref{thm:HK}]
Write $M=\del\Om$ and define
\eq{\de=\fr{\int_{M}\fr{\vt'}{H_{1}}}{\int_{M}u}-1.}
We need the following Reilly-type formula due to Qiu/Xia \cite[equ.~(16)]{QiuXia:/2015}. There holds for a arbitrary $C^{2}(\bar\Om)$-function $f$, recall \Cref{sec:notation} for notation,
\eq{\label{HK-stable-1}&\int_{\Om}\vt'(\bar\De f+(n+1)Kf)^{2}-\int_{\Om}\vt'\abs{\bar\n^{2}f+Kf \bar g}^{2}\\
	=~&\int_{M}\vt'(2\del_{\nu}f\De f+nH_{1}(\del_{\nu}f)^{2}+h(\n f,\n f)+2nKf\del_{\nu}f)\\
	&+~\int_{M}\bar\n_{\nu}\vt'(\abs{\n f}^{2}-nK f^{2}).}
Let $f$ be a solution to \eqref{HK-Diri}.
With the help of 
\eq{\bar\De\vt'=-(n+1)K\vt',}
partial integration and the divergence theorem we deduce
\eq{\int_{\Om}\vt'=\int_{M}\vt'\del_{\nu}f.}
Hence from the H\"older inequality we get
\eq{\br{\int_{\Om}\vt'}^{2}&\leq \int_{M}\fr{\vt'}{H_{1}}\int_{M}\vt' H_{1}(\del_{\nu}f)^{2}\\
			&=\int_{M}\fr{\vt'}{H_{1}}\br{\fr{1}{n+1}\int_{\Om}\vt'-\fr 1n\int_{\Om}\vt'\abs{\mr{\bar\n}^{2}f}^{2}}}
and thus
\eq{\fr{\br{\int_{\Om}\vt'}^{2}}{\fr{1}{n+1}\int_{\Om}\vt'-\fr{1}{n}\int_{\Om}\vt'\abs{\mr{\bar\n}^{2}f}^{2}}\leq \int_{M}\fr{\vt'}{H_{1}}= (1+\de)\int_{M}u=(1+\de)(n+1)\int_{\Om}\vt'.}
Rearranging gives
\eq{\int_{\Om}\vt'\abs{\mr{\bar\n}^{2}f}^{2}\leq \fr{\de n}{(1+\de)(n+1)}\int_{\Om}\vt' = \fr{n}{(1+\de)(n+1)^{2}}\int_{M}\br{\fr{\vt'}{H_{1}}-u}.}
Schauder theory \cite{GilbargTrudinger:/2001} implies that the $C^{2}$-norm of $f$ is controlled by the $C^{2,\be}$-regularity of $M$ and we note that $M$ can be covered by a controlled number of charts, see \cite[Ch.~1]{Perez:/2011}. The interior sphere condition and the gradient estimate from \autoref{quant-Hopf} gives, in conjunction with the $C^{2}$-control on $f$, a one-sided neighbourhood $\cU$ of some definite positive size, on which all the crucial quantities from \autoref{gen-Hessian-stable} are under control. Also note that $\vt'\geq \ep>0$, due to the bound on $\max_{M}r$.
Plugging this into \Cref{gen-Hessian-stable} with $p = n+1$ gives the result after estimation of $\vt'$ from below and pulling out $\abs{\bar\n^{2}f}_{0,\Om}^{p-1}$ of the integral.
\end{proof}

\subsection{Almost constant curvature functions}\label{stab-CFC}

We start with the proof of \Cref{thm:CMC}, which is related to the proof of \Cref{thm:HK}.
\begin{proof}[Proof of \Cref{thm:CMC}]
 Let $f$ be a solution to \eqref{HK-Diri} and then use $f$ in the generalised Reilly formula \eqref{HK-stable-1}. We obtain
\eq{\int_{\Om}\vt'\abs{\mr{\bar\n}^{2}f}^{2}&=\fr{n}{n+1}\int_{\Om}\vt'-n\int_{M}\vt' H_{1}(\del_{\nu}f)^{2}\\
							&=\fr{n}{(n+1)\int_{\Om}\vt'}\br{\int_{\Om}\vt'}^{2}-n\int_{M}\vt' H_{1}(\del_{\nu}f)^{2}\\
							&=\fr{n}{(n+1)\int_{\Om}\vt'}\br{\int_{M}\vt'\del_{\nu}f}^{2}-n\int_{M}\vt' H_{1}(\del_{\nu}f)^{2}\\
							&\leq\fr{n\int_{M}\vt'}{(n+1)\int_{\Om}\vt'}\int_{M}\vt'(\del_{\nu}f)^{2}-n\int_{M}\vt' H_{1}(\del_{\nu}f)^{2}\\
							&=\int_{M}\vt'(\cH - H_{1})(\del_{\nu}f)^{2},}
where
\eq{\cH = \fr{n\int_{M}\vt'}{(n+1)\int_{\Om}\vt'}.}
The $C^{2}$-norm of $f$ as well as the gradient are controlled from above and below in the same way as in the proof of \Cref{thm:HK}. Hence the result follows from \Cref{gen-Hessian-stable}.
\end{proof}

Rigidity results involving curvature functionals are often achieved by applying Newton's inequality for the normalized elementary symmetric polynomials of the principal curvatures
\eq{H_{k+1}H_{k-1}\leq H_{k}^{2},}
\cite{Newton:/1707}, also see \cite{HardyLittlewoodPolya:/1934} for the case of arbitrary $\ka\in \bbR^{n}$. If $\ka\in \Ga_{k}$, the case of equality at a point $x\in M$ occurs precisely if the point $x$ is umbilic, i.e. all principal curvatures are the same. In order to extract a stability result via the method of umbilicity, we need to control the error term in Newton's inequality and relate it to the traceless second fundamental form. The following lemma provides such control:

\begin{lemma}\label{Strict-Newton}
Let $n\geq 2$, $(\ka_{i})_{1\leq i\leq n}\in \Ga_{k}$ with
\eq{\ka_{1}\leq \dots\leq \ka_{n}.}
Let
\eq{H_{k}=\fr{1}{\binom{n}{k}}\sum_{1\leq i_{1}<\dots<i_{k}\leq n}\ka_{i_{1}}\dots\ka_{i_{k}},\q 1\leq k\leq n,}
 $H_{0}=1$ and $H_{n+1} = 0$.
Then there holds
\eq{H_{k}^{2}-H_{k+1}H_{k-1}\geq c_{n}\abs{\mr{A}}^{2}H_{k+1,n1}^{2},\q\fa 1\leq k\leq n-1,}
 where
 \eq{H_{k+1,ij} = \fr{\del^{2}H_{k+1}}{\del\ka_{i}\del\ka_{j}}.}
 Furthermore, on $\Ga_{k}$ the function $H_{k+1,n1}$ is positive. 
\end{lemma}

\pf{
From a simple calculation already done in \cite{ChenGuanLiScheuer:/2022} we obtain 
\eq{H_{k}^{2}-H_{k+1}H_{k-1}\geq c_{n}\sum_{i,j}(\ka_{i}-\ka_{j})^{2}H_{k+1,ij}^{2}\geq c_{n}\abs{\mr{A}}^{2}H_{k+1,n1}^{2}.}
 Denote by $\si_{k}$ the elementary symmetric polynomial.
Since $\ka\in \Ga_{k}$, we have
\eq{\label{Strict-Newton-1}0<\si_{k}=\ka_{n}\si_{k,n}+\si_{k+1,n}=\ka_{n}\ka_{1}\si_{k,n1}+(\ka_{n}+\ka_{1})\si_{k+1,n1}+\si_{k+2,n1}}
and
\eq{\label{Strict-Newton-2}0<\si_{k,1} = \ka_{n}\si_{k,n1}+\si_{k+1,n1}.}

Now we suppose by contradiction, that at some $\ka\in \Ga_{k}$ we had
\eq{\si_{k+1,n1}(\ka) = 0.}
First of all, this is only possible in case $n\geq 3$ and the smallest eigenvalue is $\ka_{1}\leq 0$.
Write $\la = (\ka_{2},\dots,\ka_{n-1})$ and let $S_{m}$ be the $m$-th elementary polynomial in $\la$. Then \eqref{Strict-Newton-1} implies
\eq{0<\ka_{n}\ka_{1}S_{k-2} + S_{k}}
and from \eqref{Strict-Newton-2} we obtain $S_{k-2}>0$, which gives $k\geq 2$ and $S_{k}>0$. In case $k=n-1$ this is a contradiction. For $k\leq n-2$,
 we obtain the contradiction
\eq{0<S_{k-2}S_{k}\leq c_{n,k} S_{k-1}^{2}=0. }
}

\begin{proof}[Proof of \Cref{thm:CFC}]

For $0\leq \ell\leq k$ set
\eq{\cH=\fr{\int_{M}\vt'H_{k}}{\int_{M}\vt'H_{\ell-1}},}
where here and in the following we understand
\eq{H_{-1}=\fr{1}{H_{1}}.}
We calculate for $p>n$
\eq{\int_{M}\abs{\mr{A}}^{p}&\leq C\int_{M}\abs{\mr{A}}^{2}\\
				&\leq\fr{C}{\min_{M} H_{k+1,n1}^{2}}\int_{M} H_{k+1,n1}^{2}\abs{\mr{A}}^{2}\\
				&\leq\fr{C\max_{M}H_{k}}{\min \vt'\min H_{k+1,n1}^{2}}\int_{M}\vt'\br{H_{k}-\fr{H_{k+1}H_{k-1}}{H_{k}}}\\
				&=\fr{C\max_{M}H_{k}}{\min_{M} \vt'\min_{M} H_{k+1,n1}^{2}}\int_{M}\vt'H_{\ell-1}\br{\cH-\fr{H_{k+1}}{H_{\ell}}}\\
				&\leq C\left\|\vt'H_{\ell-1}\br{\cH-\fr{H_{k+1}}{H_{\ell}}}_{+}\right\|_{1,M}. } From \cite[Thm.~2]{RothScheuer:04/2017} and the proof of \cite[Thm.~1.1]{RothScheuer:12/2018} we obtain that $M$ is embedded. The proof is complete after invoking \cite[Thm.~1.3]{De-RosaGioffre:/2021}, which was proven for hypersurfaces of the Euclidean space with $\abs{M}=1$. But as in the proof of \autoref{gen-Hessian-stable}, it is straightforward to transfer this result to hypersurfaces of conformally flat manifolds, as long as the conformal factor and the area are under control.
\end{proof}

\subsection{Stability in non-convex Alexandroff-Fenchel inequalities}\label{stab-AF}
In order to prove \Cref{thm:AF}, we use a particular case of the inverse curvature flows studied by Gerhardt \cite{Gerhardt:/1990} and Urbas \cite{Urbas:/1990}. Their result goes as follows. Given a smooth closed and starshaped hypersurface $M\sub\bbR^{n+1}$, such that the principal curvatures satisfy $\ka\in \Ga_{k}$, there exists a time-dependent family of starshaped hypersurfaces, parametrized over a sphere,
\eq{\ti x\cn [0,\8)\x \bbS^{n}\ra \bbR^{n+1},}
such that the flow equation
\eq{\del_{t}\ti x=\fr{H_{k-1}}{H_{k}}\nu}
is satisfied and the rescaled hypersurfaces
\eq{ x= e^{-t}\ti x}
smoothly converge to a sphere. The rescaled flow evolves, up to tangential repara\-metrisation, according to the flow equation
\eq{\label{GLF}\del_{t}x=\br{\fr{H_{k-1}}{H_{k}}-u}\nu.}
Guan/Li \cite{GuanLi:08/2009} have proved the Alexandroff-Fenchel inequalities for $M$ with the help of Gerhardt's and Urbas' result, but they used a slightly different rescaling of $\ti x$. However, the same result can be proved by using \eqref{GLF}. A detailed account of the flow \eqref{GLF}, which also provides an alternative argument for Gerhardt's and Urbas' result, is given in \cite{Scheuer:/2021}. There it is also shown that \eqref{GLF} keeps the quermassintegral $W_{k}(\Om_{t})$ fixed and decreases $W_{k+1}(\Om_{t})$. Here 
\eq{\del\Om_{t}=x(t,\bbS^{n})=:M_{t}.}
Since \eqref{GLF} converges to a round sphere, the proof of the Alexandroff-Fenchel inequalities is complete. 

In this section we make this proof quantitative and prove \Cref{thm:AF}.
\begin{proof}[Proof of \Cref{thm:AF}]
We may assume $M$ is smooth. The result then follows from approximation in $C^{2}$.
We are given a starshaped hypersurface $M=\del\Om\sub \bbR^{n+1}$ with $\ka\in \Ga_{k}$, and define $\ep>0$ by
\eq{\fr{W_{k+1}(\Om)}{W_{k+1}(B)}= \br{\fr{W_{k}(\Om)}{W_{k}(B)}}^{\fr{n-k}{n-k+1}}+\ep. }
For $1\leq k\leq n+1$ define
\eq{\ti W_{k}(\Om)=\fr{W_{k}(\Om)}{W_{k}(B)}.}
Let $(M_{t})_{0\leq t<\8}$ be the solution to \eqref{GLF} with $M_{0}=M$ and suppose without loss of generality that
\eq{M_{t}\ra \del B\q \mbox{in}~C^{\8},\q t\ra \8. } Then we use the variation formulae for the $\ti W_{k}$,
\eq{\fr{d}{dt}\ti W_{k+1}=c\int_{M_{t}}H_{k+1}\br{\fr{H_{k-1}}{H_{k}}-u},}
see for example \cite[Lemma~5.2]{Scheuer:/2021}. Here $c=c_{n,k}$ is a constant. We obtain		
\eq{\label{pf:AF-1}c\int_{0}^{\8}\int_{M_{t}}\br{\fr{H_{k+1}H_{k-1}}{H_{k}}-H_{k}}&=c\int_{0}^{\8}\int_{M_{t}}\br{\fr{H_{k+1}H_{k-1}}{H_{k}}-uH_{k+1}}\\
					&=\int_{0}^{\8}\fr{d}{dt}\ti W_{k+1}(\Om_{t})~dt\\
					&=1-\ti W_{k+1}(\Om)\\
					&=\ti{W}_{k}(\Om)^{\fr{n-k}{n-k+1}}-\ti W_{k+1}(\Om)\\
					&=-\ep.
}
Denote by $r$ the function on the unit sphere, by which $M$ is graphically para\-metrised. Further, write
\eq{F=\fr{H_{k}}{H_{k-1}}.} From the estimates in \cite[Lemma~4.6]{Scheuer:/2021} we see that there exists a constant $C$, which depends on $\osc r(0,\cdot)$, $(\min u(0,\cdot))^{-1}$, $1/F(0,\cdot)$ and $\abs{A}_{0,M}$ such that
\eq{\abs{A}_{0,M_{t}}+\left|\tfr 1F-u\right|\leq C.}
Hence
\eq{\dist(M_{t},M)\leq Ct\q \fa t\geq 0.}
From \eqref{pf:AF-1} we get
\eq{ \min_{s\in [0,\rt{\ep}]}\int_{M_{s}}\br{H_{k}-\fr{H_{k+1}H_{k-1}}{H_{k}}}\leq c\rt{\ep}.}
Let $M^{\ep}$ be the hypersurface where the minimum on the left is attained, in particular we obtain
\eq{\label{pf:AF-2}\dist(M^{\ep},M)\leq C\rt{\ep}\q\mbox{and}\q \int_{M^{\ep}}\br{H_{k}-\fr{H_{k+1}H_{k-1}}{H_{k}}}\leq C\rt{\ep}.}
Then, using \Cref{Strict-Newton},
\eq{\int_{M^{\ep}}\abs{\mr{A}}^{n+1}\leq C\int_{M^{\ep}}\abs{\mr{A}}^{2}&\leq \fr{C\max_{M^{\ep}}H_{k}}{\min_{M^{\ep}}H_{k+1,n1}}\int_{M^{\ep}}\fr{ H_{k+1,n1}^{2}\abs{\mr{A}}^{2}}{H_{k}}\leq C\rt{\ep}.\\
}
The proof is complete after combining \eqref{pf:AF-2} with \cite[Thm.~1.3]{De-RosaGioffre:/2021}.
\end{proof}

\subsection{Stability in Serrin's problem with nonlinear source}

The proof of \Cref{thm:Serrin} proceeds in the spirit of the integral approaches involving so-called {\it{$P$-functions}}. The crucial ingredient is the following integral identity, which is a refined version of the ones in 
\cite{BrandoliniNitschSalaniTrombetti:09/2008,MagnaniniPoggesi:/2020}. 

\begin{lemma}\label{lemma:Serrin}
Let $\Om\Subset \bbM_{0}$ be an open set with $C^{2}$-boundary and suppose $\phi$ is a $C^{2}$-function on an interval of the real line. Let $f\in C^{2}(\Om)\cap C^{1}(\bar\Om)$ satisfy
\eq{\label{eq:Serrin}\bar\De f & = \phi(f)\q&&\mbox{in}~\Om\\
			f &= 0 \q&&\mbox{on}~M=\del\Om.}
Then we have the equation
\eq{\int_{\Om}(-f)\abs{\mr{\bar\n}^{2}f}^{2} 	&=\fr{1}{2}\int_{M}(\del_{\nu}f-\del_{\nu}q)\br{\abs{\bar\n f}^{2}-R^{2}} +\fr{R^{2}}{2}\int_{\Om}(\phi - \phi(0))\\
				&\hp{=}+\int_{\Om}(\phi - \phi(0))\br{\fr{3}{2}\Phi - \fr{n}{n+1}f\phi} + \fr{\phi(0)}{2}\int_{\Om}(\Phi-f\phi),}
where
\eq{\Phi(f) = \int_{0}^{f}\phi(s)~ds,\q R = \fr{1}{\abs{M}}\int_{M}\del_{\nu}f, \q q(x) = \fr{\phi(0)}{2(n+1)}\abs{x}^{2}.}				
				\end{lemma}

\pf{
We write $\phi = \phi\circ f$.
Accounting for the nonlinearity, we define the $P$-function as
\eq{P = \fr{1}{2}\abs{\bar\n f}^{2}-\fr{\phi}{n+1}f.}
Then
\eq{\bar\De P &=d\phi(\bar\n f)+\abs{\bar\n^{2}f}^{2}-\fr{1}{n+1}f\bar\De\phi  - \fr{2}{n+1}d\phi(\bar\n f) - \fr{1}{n+1}(\bar\De f)^{2}\\
			&=\fr{n-1}{n+1}d\phi(\bar\n f) + \abs{\mr{\bar\n}^{2}f}^{2} - \fr{1}{n+1}f\bar\De\phi. }
Thus, using integration by parts,
\eq{\int_{\Om}(-f)\abs{\mr{\bar\n}^{2}f}^{2} &= \int_{\Om}(-f)\bar\De P + \fr{n-1}{n+1}\int_{\Om}fd\phi(\bar\n f) - \fr{1}{n+1}\int_{\Om}f^{2}\bar\De \phi\\
				&=\fr{1}{2}\int_{M}\del_{\nu}f\abs{\bar\n f}^{2}-\int_{\Om}\phi P + \int_{\Om}fd\phi(\bar\n f)\\
				&=\fr{1}{2}\int_{M}(\del_{\nu}f-\del_{\nu}q)\br{\abs{\bar\n f}^{2}-R^{2}} -\int_{\Om}\phi P + \int_{\Om}fd\phi(\bar\n f)\\
				&\hp{=}+ \fr{1}{2}\int_{\Om}\abs{\bar\n f}^{2}\del_{\nu}q + \fr{R^{2}}{2}\int_{\Om}(\phi - \phi(0)),}
where
\eq{q = \fr{\phi(0)}{2(n+1)}\abs{x}^{2}.}
From the Pohozaev identity, see \cite[p.~156]{Struwe:/1990}, we obtain

\eq{\int_{\Om}(-f)\abs{\mr{\bar\n}^{2}f}^{2} 	&=\fr{1}{2}\int_{M}(\del_{\nu}f-\del_{\nu}q)\br{\abs{\bar\n f}^{2}-R^{2}} -\int_{\Om}\phi P + \int_{\Om}fd\phi(\bar\n f)\\
				&\hp{=}-\phi(0)\int_{\Om}\Phi  + \fr{n-1}{n+1}\fr{\phi(0)}{2}\int_{\Om}f\phi+ \fr{R^{2}}{2}\int_{\Om}(\phi - \phi(0))\\
				&=\fr{1}{2}\int_{M}(\del_{\nu}f-\del_{\nu}q)\br{\abs{\bar\n f}^{2}-R^{2}} +\fr{3}{2}\int_{\Om}\Phi\phi- \fr{n}{n+1}\int_{\Om}f\phi^{2}\\
				&\hp{=}-\phi(0)\int_{\Om}\Phi  + \fr{n-1}{n+1}\fr{\phi(0)}{2}\int_{\Om}f\phi+ \fr{R^{2}}{2}\int_{\Om}(\phi - \phi(0)).}			
Further rearranging gives
\eq{\int_{\Om}(-f)\abs{\mr{\bar\n}^{2}f}^{2} & =\fr{1}{2}\int_{M}(\del_{\nu}f-\del_{\nu}q)\br{\abs{\bar\n f}^{2}-R^{2}} +\fr{3}{2}\int_{\Om}\Phi\br{\phi - \fr 23\phi(0)}\\
				&\hp{=}- \fr{n}{n+1}\int_{\Om}f\phi\br{\phi - \fr{n-1}{2n}\phi(0)} +\fr{R^{2}}{2}\int_{\Om}(\phi - \phi(0))\\ 
				&=\fr{1}{2}\int_{M}(\del_{\nu}f-\del_{\nu}q)\br{\abs{\bar\n f}^{2}-R^{2}} +\fr{R^{2}}{2}\int_{\Om}(\phi - \phi(0))\\
				&\hp{=}+\int_{\Om}(\phi - \phi(0))\br{\fr{3}{2}\Phi - \fr{n}{n+1}f\phi} + \fr{\phi(0)}{2}\int_{\Om}(\Phi-f\phi).}

}

\begin{proof}[Proof of \Cref{thm:Serrin}]
Let $f\in C^{2}(\bar\Om)$ be a solution of \eqref{eq:Serrin}. As $\phi>0$, \Cref{quant-Hopf} gives a lower bound on $\min_{\del\Om}\abs{\bar\n f}$ in terms of $\min \phi$ and $\rho$. 
Define
\eq{\ep = (\|\abs{\bar\n f}-R\|_{1,M}+\|\phi-\phi(0)\|_{1,\Om}+\|\Phi-f\phi\|_{1,\Om})^{\fr 12}.}
 Note that $f<0$ in $\Om$ due to the maximum principle. Let $\cU\sub\Om$ be a one-sided neighbourhood of $M$, such that 
 \eq{\min(\min_{\bar\cU}\abs{\bar\n f},\max_{\bar\cU}(-f))\geq C^{-1}>0,}
 which can be achieved due to the lower gradient bound on $M$, the $C^{2}$-bound and the interior sphere condition.  The level sets
\eq{M_{-t} = \{f=-t\},\q 0<t<\max_{\bar\cU}(-f),}
are connected, non-empty and a boundary component of 
\eq{\cU_{-t} = \{f<-t\}\cap \cU.}
We apply \Cref{gen-Hessian-stable} to the function $f$ as follows:
We estimate with \autoref{lemma:Serrin},
\eq{\int_{\cU_{-\ep}}\abs{\mr{\bar\n}^{2}f}^{n+1}\leq \fr{C}{\ep}\int_{\cU_{-\ep}}(-f)\abs{\mr{\bar\n}^{2}f}^{2}\leq C\ep. }
If $\ep$ is small, we obtain a constant $C$, such that for some sphere $S$,
\eq{\dist(M_{-\ep},S)\leq C\br{\int_{\cU_{-\ep}}\abs{\mr{\bar\n}^{2}f}^{n+1}}^{\fr{1}{n+2}}.}
Invoking that $f<-\ep$ in $\cU_{-\ep}$ and \Cref{lemma:Serrin}, we estimate further:
\eq{\dist(M_{-\ep},S)&\leq C\ep^{-\fr{1}{n+2}}\br{\int_{\cU_{-\ep}}(-f)\abs{\mr{\bar\n}^{2}f}^{2}}^{\fr{1}{n+2}}\\
		&\leq C\ep^{-\fr{1}{n+2}}\br{\int_{\Om}(-f)\abs{\mr{\bar\n}^{2}f}^{2}}^{\fr{1}{n+2}}\\
		&\leq C\ep^{-\fr{1}{n+2}}\br{\|\abs{\bar\n f}-R\|_{1,M}+\|\phi-\phi(0)\|_{1,\Om}+\|\Phi-f\phi\|_{1,\Om}}^{\fr{1}{n+2}}\\
		&= C\ep^{\fr{1}{n+2}}.}
As in the proof of \Cref{gen-Hessian-stable}, in particular recall \eqref{pf:gen-Hessian-stable-8}, we obtain
\eq{\dist(\del\Om,M_{-\ep})\leq C\ep.} 
Hence we complete the proof by
\eq{\dist(\del\Om,S)&\leq \dist(\del\Om,M_{-\ep})+\dist(M_{-\ep},S)\leq C\ep+C\ep^{\fr{1}{n+2}}.}
\end{proof}

\subsection{Stability in a fourth order Steklov eigenvalue problem}

The proof of \Cref{thm:Steklov} makes the rigidity case from \cite{WangXia:10/2009} quantitative.

\begin{proof}[Proof of \Cref{thm:Steklov}]
Let $w$ be an eigenfunction to the proposed equation, i.e.
\eq{\bar\De^{2}w & = 0\q\mbox{in}~\Om\\
		w=\bar\De w - \mu_{1}\del_{\nu}w&=0\q\mbox{on}~M=\del\Om.}
It is known that $w>0$ in $\Om$ as well as $\del_{\nu}w<0$ on $M$, \cite[Thm.~1]{BerchioGazzolaMitidieri:10/2006}.
Then
\eq{\mu_{1} = \fr{\int_{\Om}(\bar\De w)^{2}}{\int_{M}\br{\del_{\nu}w}^{2}}}
and plugging this into \eqref{HK-stable-1} with $K=0$, we get
\eq{\fr{n}{n+1}\mu_{1}\int_{M}(\del_{\nu}w)^{2}=\fr{n}{n+1}\int_{\Om}\br{\bar\De w}^{2} =\int_{M}nH_{1}(\del_{\nu}w)^{2}+ \int_{\Om}\abs{\mr{\bar\n}^{2}w}^{2}.}
This implies
\eq{(n+1)\int_{\Om}\abs{\mr{\bar\n}^{2}w}^{2}= n\int_{M}(\mu_{1}-(n+1)H_{1})(\del_{\nu} w)^{2}}
and the proof can be completed as in previous cases from \autoref{gen-Hessian-stable}.
\end{proof}

\section*{Acknowledgments}
This work was made possible through a research scholarship JS received from the DFG and which was carried out at Columbia University in New York. JS would like to thank the DFG, Columbia University and especially Prof.~Simon Brendle for their support.

%
\providecommand{\bysame}{\leavevmode\hbox to3em{\hrulefill}\thinspace}
\providecommand{\MR}{\relax\ifhmode\unskip\space\fi MR }
\providecommand{\MRhref}[2]{%
  \href{http://www.ams.org/mathscinet-getitem?mr=#1}{#2}
}
\providecommand{\href}[2]{#2}

\end{document}